\documentclass[12pt, reqno]{amsart}
\usepackage{amscd,amsthm,amsfonts,amssymb,amsmath,euscript}
\newcounter{statements}
\makeatletter\@addtoreset{statements}{subsection} \makeatother

\newtheorem{theorem}[statements]{Theorem}

\newtheorem{proposition}[statements]{Proposition}
\newtheorem{propositionn}[subsection]{Proposition}
\newtheorem{lemma}[statements]{Lemma}

\newtheorem{corollary}[statements]{Corollary}
\newtheorem{conjecture}[statements]{Conjecture}

\newtheorem{example-theorem}[statements]{Example-theorem}
\newtheorem{example-lemma}[statements]{Example-lemma}

\newtheorem{definition}[statements]{Definition}
\newtheorem{definitionn}[subsection]{Definition}
\theoremstyle{definition}

\newtheorem{examplee}[subsection]{Example}
\newtheorem{remark}[statements]{Remark}
\newtheorem{remarkk}[subsection]{Remark}
\theoremstyle{remark}

\newcommand{\QQ}{{\mathbb Q}}
\newcommand{\ZZ}{{\mathbb Z}}
\newcommand{\NN}{{\mathbb N}}

\newcommand{\PP}{{\mathbb P}}
\newcommand{\CC}{{\mathbb C}}
\newcommand{\Aff}{{\mathbb A}}

\newcommand{\itc}[1]{\textup{#1}}
\newcommand{\ind}{\mathrm{ind}\,}
\newcommand{\virt}{\mathrm{virt}}
\newcommand{\tit}{Gromov--Witten invariants of Fano threefolds
\\
of genera $6$ and $8$}

\begin{document}

\begin{title}
\tit
\end{title}

\begin{abstract}
    The aim of this paper is to prove Golyshev's conjecture in the cases
of Fano threefolds $V_{10}$ and $V_{14}$. This conjecture states
modularity of $D3$ equations for smooth Fano threefolds with
Picard group $\ZZ$. More precisely, we find counting matrices of
prime two-pointed Gromov--Witten invariants for them. For this we
use the method that lets us find Gromov--Witten invariants of
complete intersections in varieties whose invariants are
(partially) known.

\end{abstract}

\author{Victor Przyjalkowski}

\thanks{The work was partially supported by RFFI grants
$04-01-00613$ and $05-01-00353$, and grant NSh$-489.2003.1$.}

\address{Steklov Institute of Mathematics, 8 Gubkin street, Moscow 117966, Russia}

\email{victorprz@mi.ras.ru}

\maketitle

In this work we use the following method to find two-pointed
Gromov--Witten invariants for $V_{10}$ ð $V_{14}$. According to
S.\,Mukai, these varieties are complete intersections in
grassmannians (Theorem~\ref{Mukai}). Find their generating series
for one-pointed invariants with descendants (the so called
$I$-series). By Quantum Lefschetz Theorem (see~\ref{fano lefschetz})
find $I$-series for $V_{10}$ and $V_{14}$. Finally, by divisor axiom
and topological recursion relations find polynomial expressions for
the coefficients of $I$-series in terms of two-pointed prime
invariants (Proposition~\ref{lemma:threefold_equations}). Finally,
inverse these expressions and find two-pointed invariants.

The paper is organized as follows. In $\S\,1$ we give definitions of
Gromov--Witten invariants and stacks of stable maps of genus $0$
with marked points. In $\S\,2$ we state Golyshev's conjecture. In
$\S\,3$ we consider $I$-series for grassmannians. $\S\,4$ contains
Quantum Lefschetz Theorem. The relations between one- and
multi-pointed invariants are in $\S\,5$. Finally, in $\S\,6$ we
prove the main theorem of this paper (Theorem~\ref{third}), in which
we explicitly find the specific counting matrices for threefolds
$V_{10}$ and $V_{14}$.

\medskip

{\bf Bibliography.} Gromov--Witten invariants were introduced for
counting the numbers of curves of different genera on different
varieties, which intersect given homological classes. Such numbers
are called \emph{prime} invariants. Axiomatic treatment of them was
given in~\cite{KM1}. After that these invariants were explicitly
constructed. The generalization (the so called \emph{invariants with
descendants}) was obtained in~\cite{BM1} and~\cite{Beh}.

    The idea to express the Gromov--Witten invariants of
hypersurfaces in terms of the invariants of the ambient variety (the
so called Mirror Formula), the main (up to now) tool of finding of
quantum cohomology, seems to belong to A.\,Givental (see, for
instance,~\cite{Gi}). We use it in the form given by A.\,Gathmann
(see ~\cite{Ga2}). The formula for $I$-series for grassmannians (the
Hori--Vafa conjecture) was proved by A.\,Bertram,
I.\,Ciocan-Fontanine, B.\,Kim in~\cite{BCK}. A.\,Bertram with
H.\,Kley in~\cite{BK} and Y.-P.\,Lee with R.\,Pandharipande
in~\cite{LP} express multipointed invariants in terms of one-pointed
ones.

Fano threefolds of genera $6$ and $8$ were studied in~\cite{Gu1},
\cite{Gu2}, \cite{Gu3}, \cite{Is1}, \cite{Is2},~\cite{I},~\cite{IP},
\cite{Lo1},~\cite{Lo2},~\cite{M1},~\cite{Mar}. Gromov--Witten
invariants of such varieties that count the lines are computed
in~\cite{Mar}, see also~\cite{IP}. The number of conics on the
threefold of genus $6$ that pass through a general point is found
in~\cite{Lo1}.

Gromov--Witten invariants of Fano threefolds first was studied
in~\cite{BM2}. Quantum cohomology of some threefolds was found
in~\cite{Ci} (some blow-ups of projective spaces and
quadrics),~\cite{AM},~\cite{Ba} (toric varieties),~\cite{QR} (vector
bundles over $\PP^n$), and~\cite{Bea} (some complete intersections).

Conjecture on modularity of smooth Fano threefolds with Picard Group
$\ZZ$ is stated by V.\,Golyshev in~\cite{Go1}.

\bigskip

\paragraph{\textbf{Conventions and notations}}

    Everything
    is over $\CC$.
    Throughout the paper:

    \emph{A grassmannian} $G(r,n)$ is the variety of $r$-dimensional linear subspaces in
    the $\mbox{$n$-dimensional}$ linear space.

    In the sequel \emph{the Gromov--Witten invariants} (GW-invariants) mean
the \emph{genus $0$ invariants} (i. e. the invariants that
correspond to rational curves).

The cohomology algebra $H^*(X,\QQ)$ of variety $X$ we denote as
$H^*(X)$.

Poincar$\rm \acute{e}$ dual class to the class $\gamma\in H^*(X)$ on
$X$ we denote as $\gamma^\vee$.

The subalgebra of algebraic cycles of $H^*(X)$ we denote as
$H^*_{alg}(X)$.

The subset of classes of effective curves in $H_2(X,\mathbb Z)$ we
denote by $H^+_2(X)$.

\section{Main definitions}

    Let $X$ be a smooth projective variety with Picard group $\ZZ$ such that $-K_X$ is nef.

We mention only those definitions and axioms that will be used in
the sequel.

\subsection{Moduli spaces.}
\begin{definition}
\emph{The genus} of curve $C$ is the number $h^1(\mathcal O_C)$.
\end{definition}

It is easy to see that the curve is of genus $0$ if and only if it
is a tree of rational curves.

\begin{definition} {\it
The connected curve $C$ with $n\geq 0$ marked points
$p_1,\ldots,p_n\in C$ is called \emph{prestable} if it has at most
ordinary double points as singularities and $p_1,\ldots, p_n$ are
distinct smooth points \itc{(}see~\cite{Ma}, III--2.1\itc{)}. The
map $f\colon C\rightarrow X$ of connected curve of genus $0$ with
$n$ marked points are called \emph{stable} if $C$ is prestable and
there are at least three marked or singular points on every
contracted component of $C$ \itc{(}\cite{Ma}, V--1.3.2\itc{)}. }
\end{definition}

In the other words, a stable map of connected curve is the map that
has only finite number of infinitesimal automorphisms.

\begin{definition}
\emph{The family of stable maps} \itc{(}over the scheme $S$\itc{)}
of curves of genus $0$ with $n$ marked points is the collection
$(\pi \colon \mathcal C\rightarrow S, p_1,\ldots,p_n, f\colon
\mathcal C\rightarrow X)$, where $\pi$ is the following map. It is a
smooth projective map with $n$ sections $p_1,\ldots,p_n$. Its
geometric fibers $(\mathcal C_s,p_1(s),\ldots,p_n(s))$ are prestable
curves of genus $0$ with $n$ marked points. Finally, the restriction
$f|_{\mathcal C_s}$ on each fiber is a stable map.

Two families over $S$
$$
(\pi\colon \mathcal C\rightarrow S, p_1,\ldots,p_n, f),\ \
(\pi^\prime\colon \mathcal C^\prime\rightarrow S,
p_1^\prime,\ldots,p_n^\prime, f^\prime)
$$
are called \emph{isomorphic} if there is an isomorphism $\tau\colon
\mathcal C\rightarrow \mathcal C^\prime$ such that
$\pi=\pi^\prime\circ \tau$, $p_i^\prime=\tau\circ p_i$,
$f=f^\prime\circ\tau$.
\end{definition}

Let $\beta\in H_2^+(X)$. Consider the following (contravariant)
functor $ \overline{\mathcal M}_n(X, \beta)$ from the category of
(complex algebraic) schemes to the category of sets. Let
$\overline{\mathcal M}_n(X, \beta)(S)$ be the set of isomorphism
classes of families of stable maps of genus $0$ curves with $n$
marked points $(\pi\colon \mathcal C\rightarrow S, p_1,\ldots, p_n,
f)$ such that $f_*([\mathcal C_s])=\beta$, where $[C_s]$ is the
fundamental class of $C_s$.

\begin{definition}
{\it \emph{The moduli space of stable maps} of genus $0$ curves of
class $\beta\in H^+_2(X)$ with $n$ marked points to $X$ is the
Deligne--Mumford stack \itc{(}see~\cite{Ma}, V--5.5\itc{)} which is
the coarse moduli space that represents $\overline{\mathcal M}_n(X,
\beta)$. This space is denoted by $\bar{M}_n(X,\beta)$. }
\end{definition}

The stack $\bar{M}_n(X,\beta)$ is compact and smooth, that locally
it is a quotient of a smooth variety by a finite group. Hence, we
can consider an intersection theory on it (see~\cite{Vi}). In the
general case $\bar{M}_n(X,\beta)$ has a ``wrong'' dimension, so it
is equipped with \emph{the virtual fundamental class}
$[\bar{M}_n(X,\beta)]^\virt$ of virtual dimension $\mathrm{vdim}\,
\bar{M}_n(X,\beta)=\dim X - \deg_{K_X} \beta+n-3$. Let $X$ be a
\emph{convex} variety (i. e. for any map $\mu\colon \PP^1\rightarrow
X$ the equality $H^1(\PP^1,\mu^*(TX))=0$ holds). Then
$\bar{M}_n(X,\beta)$ is projective normal variety of pure dimension
$\dim X - \deg_{K_X} \beta+n-3$. Thus, the virtual fundamental class
is the usual one. More in~\cite{Ma}, VI--1.1.

\subsection{Gromov--Witten invariants.}

Consider \emph{the evaluation maps}
$ev_i:\bar{M}_n(X,\beta)\rightarrow X$, given by
$ev_i(C;p_1,\ldots,p_n,f)=f(p_i)$. Let
$\pi_{n+1}:\bar{M}_{n+1}(X,\beta)\rightarrow \bar{M}_n(X,\beta)$ be
\emph{the forgetful map} at the point $p_{n+1}$ (which forget this
point and contract unstable component after it); consider the
sections $\sigma_i:\bar{M}_n(X,\beta)\rightarrow
\bar{M}_{n+1}(X,\beta)$, which coincide with the points $p_i$. The
image of a curve $(C;p_1,\ldots,p_n,f)$ under $\sigma_i$ is a curve
$(C';p_1,\ldots,p_{n+1},f')$. Here $C'=C\bigcup C_0$, $C_0\backsimeq
\PP^1$, $C_0$ and $C$ intersect at the non-marked point $p_i$ on
$C'$, and $p_{n+1}$ and $p_i$ lie on
 $C_0$. The map $f'$ contracts $C_0$ to the point and $f'|_C=f$.

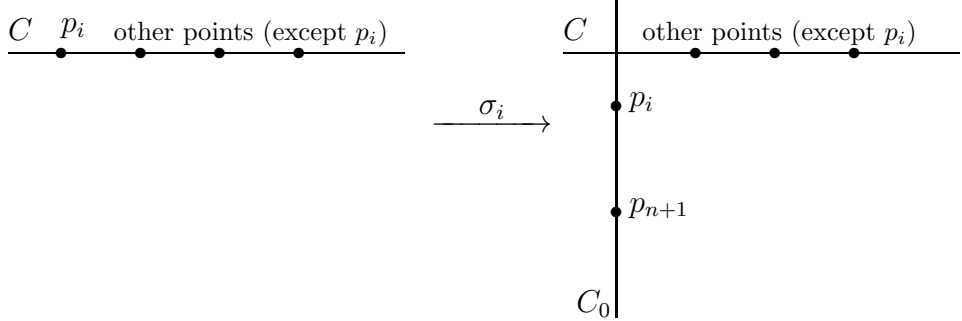
\begin{picture}(200,130)(-40,0)
\put(0,110){\line(1,0){150}} \put(0,115){$C$}
\put(50,110){\circle*{4}} \put(40,115){\footnotesize other points
(except $p_i$)} \put(20,110){\circle*{4}} \put(20,118){$p_i$}
\put(80,110){\circle*{4}} \put(110,110){\circle*{4}}

\put(178,87){$\sigma_i$} \put(160,80){$\xrightarrow{\ \ \ \ \large
\ \ \ \ \ }$}

\put(215,12){$C_0$}\put(210,110){\line(1,0){150}}
\put(210,115){$C$}\put(230,10){\line(0,1){120}}
\put(260,110){\circle*{4}} \put(240,115){\footnotesize other
points (except $p_i$)} \put(290,110){\circle*{4}}
\put(320,110){\circle*{4}} \put(230,90){\circle*{4}}
\put(235,90){$p_i$} \put(230,50){\circle*{4}}
\put(235,50){$p_{n+1}$}
\end{picture}

Consider the sheaf $L_i$ given by
$L_i=\sigma_i^*\omega_{\pi_{n+1}}$, where $\omega_{\pi_{n+1}}$ is
a relative dualizing sheaf of $\pi_{n+1}$. Its fiber over the
point $(C;p_1,\ldots,p_n,f)$ is $T^*_{p_i}C$. Put \emph{cotangent
line class} $\psi_i=c_1(L_i)\in H^2(\bar{M}_n(X,\beta))$.

\begin{definition}[\cite{Ma}, VI--2.1]  {\it
Consider $\gamma_1,\ldots, \gamma_n\in H^*(X)$ and let
$d_1,\ldots,d_n$ be non-negative integers. Then \emph{the
Gromov--Witten invariant (correlator)} is the number given by
\begin{equation*}
\langle\tau_{d_1} \gamma_1,\ldots,
\tau_{d_n}\gamma_n\rangle_\beta=
ev_1^*\gamma_1\cdot\psi_1^{d_1}\cdot\ldots \cdot
ev_n^*\gamma_n\cdot \psi_n^{d_n}\cdot [\bar{M}_n(X,\beta)]^\virt
\end{equation*}
if $\sum \mathrm{codim}\, \gamma_i+\sum d_i=\mathrm{vdim}\,
\bar{M}_n(X,\beta)$ and 0 otherwise. The number $\sum d_i$ is called
\emph{the degree} of invariant with respect to the descendants. The
invariants of degree $0$ are called \emph{prime} \itc{(}and symbols
$\tau_0$ are omitted\itc{)}. }
\end{definition}

After this definition \emph{axioms}, which define Gromov--Witten
invariants, become their \emph{properties}. We will refer to two
of them.

\begin{itemize}
  \item Divisor axiom (\cite{Ma}, V--5.4).

  Let $\gamma_0\in H^2(X)$. Then
\begin{multline}
\label{axiom:divisor} \langle\gamma_0, \tau_{d_1} \gamma_1,\ldots,
\tau_{d_n}\gamma_n\rangle_\beta=
(\gamma_0\cdot\beta)\langle\tau_{d_1} \gamma_1,\ldots, \tau_{d_n}\gamma_n\rangle_\beta+\\
\sum_{k,d_k\ge 1} \langle\tau_{d_1} \gamma_1,\ldots,
\tau_{d_k-1}(\gamma_0\cdot
\gamma_k),\ldots,\tau_{d_n}\gamma_n\rangle_\beta.
\end{multline}
The only exception is the case $\beta=0$ and $n=2$. In this case
$\langle\gamma_0,\gamma_1,\gamma_2\rangle_0=\gamma_0\cdot\gamma_1\cdot\gamma_2\cdot
[X]$.
  \item Fundamental class axiom (\cite{Ma}, V--5.1).

Let $\mathbf 1$ be the fundamental class of $X$. Then
\begin{equation}
\langle \mathbf{1}, \tau_{d_1} \gamma_1,\ldots,
\tau_{d_n}\gamma_n\rangle_\beta= \sum_{k,d_k\ge
1}\langle\tau_{d_1} \gamma_1,\ldots,
\tau_{d_k-1}\gamma_k,\ldots,\tau_{d_n}\gamma_n\rangle_\beta.
\end{equation}
The only exception is the case $\beta=0$ and $n=2$. In this case
$\langle
\mathbf{1},\gamma_1,\gamma_2\rangle_0=\gamma_1\cdot\gamma_2\cdot
[X]$.
\end{itemize}

\begin{remark}
    In this paper Gromov--Witten invariants are considered
for varieties of Picard rank $1$. In this case the class of the
curve $\beta\in H^+_2(X)$ is determined by its (anticanonical)
degree $d$, and we will often write $\langle\ldots\rangle_d$ instead
of $\langle\ldots\rangle_\beta$. For the case of greater Picard rank
one should use \emph{multidegree}.

\end{remark}

\subsection{$I$-series~(\cite{Ga2})}

Let $R\in H^*(X)$ be Poincar$\rm \acute{e}$ self-dual subalgebra
(that is, for any class $\gamma\in R$ we have $\gamma^\vee\in R$).
Let $1=\gamma_0,\ldots,\gamma_N$ be the basis in $R$. Let $\beta\in
H^+_2(X)$ be an effective curve of degree $d$ (with respect to the
positive Picard group generator). Put
\begin{gather*}
I^X_{d,R}=I^X_{\beta,R}=\sum_{i,j}
\langle\tau_i \gamma_j\rangle_\beta\gamma^\vee_j.\\
\end{gather*}

\begin{definition}[\cite{Ga2}] {\it
\emph{$I$-series} $I^X_R$ is given by the following.
\begin{gather*}
I^X_R=\sum_{d\geq 0} I_{d,R}^X\cdot q^d.
\end{gather*}         }
\end{definition}

We denote $I$-series $I^X_{R}$ for $R=H^*(X)$ by $I^X$; for
$R=H^*_{alg}(X)$ by $I^X_{alg}$. Let $Y\subset X$ be a complete
intersection and  $\pi\colon H^*(X)\rightarrow H^*(Y)$ be a natural
restriction homomorphism. Then for $R=\pi(H^*(X))$ the series
$I^X_R$ we denote by $I^X_{rest}$.

\begin{definition}
Let $I\in H^*(X)\otimes \CC[[q]]$. Put
$$
I=\sum_{0\leq i\leq N} \gamma_i \otimes I^{(i)}.
$$
The term $I^{(0)}$ we denote by $I_{H^0}$.
\end{definition}

\section{Counting matrices for Fano threefolds}

\subsection{Counting matrices}
\begin{definition}[\cite{Go1}]
Let $X$ be a smooth threefold Fano variety such that $\mathrm{Pic}\,
(X)\cong \ZZ$ and $K=-K_X$. \emph{Counting matrix} is the matrix of
its GW-invariants, namely the following matrix $A\in
\mathrm{Mat}(4\times 4)$.
$$
    \left[\begin{array}{cccc}
    a_{00} & a_{01} & a_{02} & a_{03} \\
    1      & a_{11} & a_{12} & a_{13} \\
    0      & 1      & a_{22} & a_{23} \\
    0      & 0      & 1      & a_{33} \\
    \end{array}\right].
$$
Numeration of rows and columns starts from $0$ and the elements are
given by
$$
a_{ij}=\frac{\langle K^{3-i}, K^j, K\rangle_{j-i+1}}{\deg
X}=\frac{j-i+1}{\deg X}\cdot\langle K^{3-i}, K^j \rangle_{j-i+1}$$
\itc{(}the degree is taken with respect to the anticanonical
class\itc{)}.
\end{definition}

It is easy to see that the matrix $A$ is symmetric with respect to
the secondary diagonal: $a_{ij}=a_{3-j,3-i}$. By definition,
$a_{ij}=0$ if $j-i+1<0$. If $j-i+1=0$, then $a_{ij}=1$, because it
is just a number of intersection points of $K^{3-i}$, $K^j$, and
$K$, which is $\deg X$; $a_{00}=a_{33}=0$. For the other
coefficients $a_{ij}$'s are ``expected'' numbers of rational curves
of degree $j-i+1$ passing through $K^{3-i}$ and $K^j$, multiplied by
$\frac{j-i+1}{\deg X}$. The only exception is the following: by
divisor axiom
$$
a_{01}=2\cdot (2\cdot ind\,(X)\cdot\mathrm{[\emph{the number of
conics passing through the general point}])}.
$$

\subsection{Golyshev's conjecture}
It is more convenient to use \emph{a family} $A^\lambda=A+\lambda
E$, where $E$ is the identity matrix. Thus, the element of the
family $A^\lambda$ is given by \emph{six} parameters: five different
GW--invariants $a_{ij}$ and $\lambda$.

Consider the one-dimensional torus $\mathbb G_m= \mathrm{Spec}\, \CC
[t, t^{-1}]$ and the differential operator
$D=t\frac{\partial}{\partial t}$. Construct the family of matrices
$M^\lambda$ in the following way. Put its elements $m_{kl}^\lambda$
as follows:
$$
    m_{kl}^\lambda=
  \begin{cases}
    0, & \text{if $k>l+1$}, \\
    1, & \text{if $k=l+1$}, \\
    a_{kl}^\lambda\cdot (Dt)^{l-k+1}, & \text{if $k<l+1$}.
  \end{cases}
$$
Now consider the family of differential operators
$$
    \widetilde{L}^\lambda=\mathrm{det}_{\mathrm{right}}(DE-M^\lambda),
$$
where $\mathrm{det}_{\mathrm{right}}$ means ``right determinant'',
i. e. the determinant, which is calculated with respect to \emph{the
rightmost} column; all minors are calculated in the same way.
Dividing $\widetilde{L}^\lambda$ on the left by $D$, we get the
family of operators $L^\lambda$, so
$\widetilde{L}^\lambda=DL^\lambda$.
\begin{definition}[\cite{Go1}, 1.8]
    {\it The equation of the family $L^\lambda[\Phi(t)]=0$ is called \emph{counting equation
    $D3$}.  }
\end{definition}

\begin{conjecture}[V.\,Golyshev,~\cite{Go1}]
\label{conjecture} The solution of the $D3$ equation for smooth Fano
threefold is modular. More precisely, let $X$ be such variety, $i_X$
be its index, and $N=\frac{\deg X}{2i_X^2}$. Then in the family of
counting equations for $X$ there is one, $L^{\lambda_X}[\Phi(t)]=0$,
whose solution is an Eisenstein series of weight $2$ on $X_0(N)$.
\end{conjecture}

Based on this conjecture, V.\,Golyshev in~\cite{Go1} gives a list of
predictions for counting matrices of Fano threefolds. So we should
find all of the counting matrices for Fano varieties to check it.

\section{$I$-series and grassmannians}

\subsection{Mukai Theorem}
\begin{theorem}[\cite{M1}]
\label{Mukai}
    A smooth Fano threefold $V_{10}$ of genus $6$ \itc{(}and anticanonical degree
$10$\itc{)} is a section of grassmannian $G(2,5)$ by a linear
subspace of codimension 2 and a quadric in the Pl\"{u}cker
embedding.

    A smooth Fano threefold $V_{14}$ of genus $8$ \itc{(}and anticanonical degree
$14$\itc{)} is a section of grassmannian $G(2,6)$ by a linear
subspace of codimension 5 in the Pl\"{u}cker embedding.
\end{theorem}

    This theorem is a particular case of the general Mukai Theorem,
which describes all smooth Fano threefolds with Picard number 1 as
sections of certain sheaves on grassmannians.

\subsection{$I$-series for grassmannians}

\begin{theorem}[The Hori--Vafa conjecture (\cite{HV}, Appendix A), proof in~\cite{BCK}]
\label{theorem:HV} Let $x_1,\ldots x_r$ be the Chern roots of the
dual to the tautological subbundle $S^*$ on $G=G(r,n)$ and $r>1$.
Then
\begin{equation*}
\label{equation:HV} I^G=\sum_{d\ge 0} (-1)^{(r-1)d}
\sum_{d_1+\ldots+d_r=d}\frac{\prod_{1\le i< j\le
r}(x_i+d_i-x_j-d_j)} {\prod_{1\le i< j\le r}(x_i-x_j)
\prod_{i=1}^r\prod_{l=1}^{d_i}(x_i+l)^n} q^d.
\end{equation*}
In the case $r=1$ \itc{(}i. e. for projective space\itc{)}
\begin{equation*}
    I^{\PP^{n-1}}=\sum_{d\ge 0}\prod_{i=1}^{d} \frac{q^d}{(H+i)^n},
\end{equation*}
where $H$ is Poincar$\acute{e}$ dual to the class of hyperplane
section.
\end{theorem}

\begin{corollary}[\cite{BCK}, Proposition $3.5$]
\label{corollary:free term} The constant term of $I^G$ for
$G=G(2,n)$ is
$$
\sum_{d\geq 0}\frac{q^d}{(d!)^n}\frac{(-1)^d}{2}\sum_{m=0}^{d}\left(
\begin{array}{c}
d\\
m
\end{array}\right)^n
\left( n(d-2m)(\gamma (m)-\gamma (d-m))+2 \right),
$$
where $\gamma (m)=\sum_{j=1}^m\frac{1}{j}$ and $\gamma (0)=0$.
\end{corollary}

\section{Quantum Lefschetz Theorem}

The variety $X$ is smooth projective and with Picard group $\ZZ$ as
before.

\subsection{$I$-series of hypersurface}

\begin{lemma}[\cite{Ga1}, Lemma 5.5 or proof of
Lemma 1 in~\cite{LP}]
\label{lemma:non-algebraic} Let $Y\subset X$ be
a complete intersection and $\varphi \colon H^*(X) \rightarrow
H^*(Y)$ be the restriction homomorphism. Let
$\widetilde{\gamma}_1\in \varphi (H^*(X))^\bot$ and
$\gamma_2,\ldots, \gamma_l\in \varphi (H^*(X))$. Then for each
$\beta \in \varphi (H^+_2(X))\subset H^+_2(Y)$ the Gromov--Witten
invariant on $Y$ of the form
$$
\langle \tau_{d_1} \widetilde{\gamma}_1, \tau_{d_2} \gamma_2, \ldots
\tau_{d_l} \gamma_l \rangle_\beta
$$
vanishes.
\end{lemma}

Thus, by divisor axiom~(\ref{axiom:divisor}) and this lemma
$I^Y=I^Y_{rest}$.

\begin{remark}
In contrast to the one-pointed invariants, the two-pointed ones,
which correspond the classes in $\varphi (H^*(X))^\bot$, can be
non-zero (see.~\cite{Bea}, Proposition 1).
\end{remark}

\subsection{Mirror Formula}

\begin{theorem}[Mirror Formula,~\cite{Ga2}, Corollary 1.13]
\label{theorem:MF}
    Let $Y\subset X$ be a hypersurface and $-K_Y\ge 0$. Then
there exist series $R(q)\in H^*(X)[[q]]$ and $S(q)\in H^*(X)[[q]]$
such that
\begin{equation*}
    \sum_\beta\prod_{i=0}^{Y\cdot\beta}(Y+i)\cdot I_\beta^X\cdot q^{Y\cdot\beta}=
    R(q)\cdot \sum_\beta I_\beta^Y\cdot \widetilde{q}^{Y\cdot\beta},
\end{equation*}
where $\widetilde{q}=q\cdot e^{S(q)}$.
\end{theorem}

A.\,Gathmann describes these series in Definition $1.11$ and Lemma
$1.12$ from~\cite{Ga2}.

 In the case of Fano varieties the Mirror Formula may
be simplified.

\begin{corollary}
\label{fano lefschetz}
    Let $I=\ind Y\in \NN$ be the index of Fano variety $Y$.
If $Y$ is complete intersection of hypersurfaces $Y_1,\ldots,Y_k$ of
degrees $d_1,\ldots,d_k$ in $X$, then
\begin{equation}
\label{MF}
   \sum_\beta\prod_{j=1}^k\prod_{i=0}^{Y_j\cdot\beta}(Y_j+i)\cdot I_\beta^X\cdot
   q^{Y\cdot\beta}=
   e^{\alpha_Y\cdot q^{\deg Y}}\cdot \sum_\beta I_\beta^Y\cdot q^{Y\cdot\beta},
\end{equation}
where $\alpha_Y=\prod d_i!\cdot I_{1,H^0}^X$ if $I=1$ and
$\alpha_Y=0$ if $I\ge 2$.
\end{corollary}

\begin{proof}
The form of series $R(q)$ and $S(q)$ involves that for Fano
varieties $R(q)=e^{\alpha_Y \cdot q^{\deg Y}}$ and $S(q)=0$. For
$I>1$ we have $R(q)=1$. To find the number $\alpha_Y$ for the other
case notice that by the dimensional argument $\langle H^{\dim Y}
\rangle_l=0$, where $H$ is a hyperplane section of $X$. Comparing
the coefficients of the formula~(\ref{MF}) at $q^{\deg Y}$ and
$H^k$, we get
$$
    \alpha_Y=\prod_{i=1}^k d_k!\cdot I_{1,H^0}^X.
$$
\end{proof}

\begin{remark}
    Conjecture~\ref{conjecture} may be generalized in the following way.
Put $\alpha_Y=\lambda_Y$. Then the equation
$L^{\lambda_X}[\Phi(t)]=0$ has an Eisenstein series of weight 2 on
$X_0(N)$ as a solution.
\end{remark}

\section{Expressions for one-pointed invariants in terms of prime
two-pointed ones}

\begin{definitionn}
Subalgebra $R\in H^*(Y)$ is called quantum self-dual if for all
$\gamma\in R$, $\mu\in R^\bot$ and $\beta \in H_2^+(Y)$ the
following holds: $\gamma^\vee\in R$ and $\langle \gamma,
\mu\rangle_\beta=0$.
\end{definitionn}

\begin{propositionn}
\label{proposition:simple expressions} For each $n, I\in \NN$,
$k,d\in \ZZ_{\geq 0}$ there exist polynomial $f_k^d\in
\QQ[a_{i,j}]$, $0 \leq i,j\leq n$, $j-i+1\leq d$ such that the
following holds. Consider Fano variety $Y$ of dimension $n$ and
index $I$ such that subalgebra $R\in H^*(Y)$ generated by $-K_Y=IH$
is quantum self-dual. Then
$$
\langle\tau_k H^{d+n-2-k}\rangle_d=f_k^d(a_{ij}),
$$
where $a_{i,j}=\langle H^{n-i}, H^\beta, H \rangle_{j-i+1}/\deg Y=
\frac{j-i+1}{\deg Y}\langle H^{n-i}, H^j \rangle_{j-i+1}$.
\end{propositionn}

\begin{proof}
The strategy of finding $f_k^d$ is the following: express given
one-pointed invariant in terms of three-pointed ones (with
descendants) by divisor axiom. Then, by topological recursion
relations, express these three-pointed invariants in terms of
two-pointed ones.

Applying the divisor axiom for $H$
\begin{equation*}
\langle\tau_k H^{d+n-2-k}\rangle_d=1/d\cdot(\langle H,\tau_k
H^{d+n-2-k}\rangle_d- \langle\tau_{k-1} H^{d+n-1-k}\rangle_d)
\end{equation*}
we get
\begin{equation*}
\langle \tau_k H^{d+n-2-k}\rangle_d= \frac{1}{d}
\sum_{i=0}^k\frac{(-1)^i}{d^i} \langle H,\tau_{k-i}
H^{d+i+n-2-k}\rangle_d.
\end{equation*}

Now, by topological recursion (see, for instance,~\cite{KM2},
Corollary $1.3$)
\begin{multline*}
\langle H^a,\tau_k H^{d+n-1-a-k}\rangle_d= 1/d(\sum_{d_1\le d}
\langle H^{d_1-d+1+a}, \tau_{k-1} H^{d+n-1-a-k}\rangle_{d_1}\cdot
a_{d_1-d+a+1, a} -\\
\langle H^a,\tau_{k-1} H^{d+n-a-k}\rangle_d).
\end{multline*}

This formula may be simplified further. If we express the last
summand recursively on the right, we get
\begin{multline*}
\langle H^a, \tau_k H^{d+n-1-a-k}\rangle_d= \sum_{
\parbox[c]{1,2cm}{\scriptsize
$i=1..k,$ $d_1\le d$}} \frac{(-1)^{i+1}}{d^i}a_{d_1-d+1+a,a}
\langle H^{d_1-d+1+a}, \tau_{k-i} H^{d+n-2-a-k+i}\rangle_{d_1}\\
+\frac{(-1)^k}{d^k}\langle H^a, H^{d+n-a}\rangle_d.
\end{multline*}

\end{proof}

\begin{remarkk}
\label{remark:general_expressions}

We make the assumption that $\mathrm{Pic}\, Y= \ZZ$ and $R$ is
generated by $H$ for simplicity, as we use this in the following.
One can find the similar expressions in the general case under the
only assumption $R\cap H^2(Y)\neq\emptyset$, using the same
arguments as in the proof of Proposition~\ref{proposition:simple
expressions}.
\end{remarkk}

\begin{examplee}
\label{example:expressions} Consider a smooth Fano variety $Y$ of
index $1$ and dimension $3$, whose algebraic cohomology is generated
by a generator of Picard group $H$. Then we have:
\begin{itemize}
  \item The constant term of $I$-series.
\begin{description}
  \item[1] $\langle H^{3}\rangle_{2}=\deg Y\cdot a_{01}/4;$
  \item[2] $\langle \tau H^3 \rangle_{3}=\deg Y\cdot (a_{11}a_{01}/18+a_{02}/27);$
  \item[3] $\langle \tau_{2} H^3\rangle_{4}=\deg Y\cdot (a_{01}^2/64+a_{11}^2a_{01}/96+7 a_{11}a_{02}/576+
  a_{01}a_{12}/128+a_{03}/256).$
\end{description}
  \item The linear term of $I$-series with respect to $H$.
\begin{description}
  \item[1] $\langle H^{2}\rangle_{1}=\deg Y\cdot a_{11};$
  \item[2] $\langle \tau H^{2} \rangle_{2}=\deg Y\cdot (a_{11}^2/4+a_{12}/8-a_{01}/4);$
  \item[3] $\langle \tau_{2} H^2 \rangle_{3}=\deg Y\cdot (5 a_{11}a_{01}/108+a_{11}^3/18+a_{11}a_{12}/12-2a_{02}/81);$
  \item[4] $\langle \tau_{3} H^2\rangle_{4}=\deg Y\cdot (13a_{11}^2a_{01}/576+17a_{11}a_{02}/1728-a_{03}/256-
  3a_{01}^2/128+a_{11}^4/96+a_{12}^2/256+a_{11}^2a_{12}/32).$
\end{description}
\end{itemize}
Thus,

$
I^Y=1+a_{11}q+(a_{01}/4+(a_{11}^2/4+a_{12}/8-a_{01}/4)H)q^2+((a_{11}a_{01}/18+a_{02}/27)+(5
a_{11}a_{01}/108+a_{11}^3/18+a_{11}a_{12}/12-2a_{02}/81)H)q^3+((a_{01}^2/64+a_{11}^2a_{01}/96+7
a_{11}a_{02}/576+
  a_{01}a_{12}/128+a_{03}/256)+(13a_{11}^2a_{01}/576+17a_{11}a_{02}/1728-a_{03}/256-
  3a_{01}^2/128+a_{11}^4/96+a_{12}^2/256+a_{11}^2a_{12}/32)H)q^4+\ldots
\ \, (\mathrm{mod}\ \, H^2). $
\end{examplee}

\section{Main theorem}

\subsection{Complete intersections in grassmannians}

\begin{theorem}
\label{third}
    The counting matrices of Fano threefolds $V_{10}$ and $V_{14}$ coincide with
the predictions in~\cite{Go1}.
\begin{description}
  \item[1] For $V_{10}$
$$
    M(V_{10})=
    \left[\begin{array}{cccc}
         0 &    156 &   3600 &  33120 \\
    1      &     10 &    380 &   3600 \\
    0      & 1      &     10 &    156 \\
    0      & 0      & 1      &      0 \\
    \end{array}\right].
$$
The shift $\alpha_{V_{10}}$ is $6$.

  \item[2] For $V_{14}$
$$
    M(V_{14})=
    \left[\begin{array}{cccc}
         0 &     64 &    924 &   5936 \\
    1      &      5 &    140 &    924 \\
    0      & 1      &      5 &     64 \\
    0      & 0      & 1      &      0 \\
    \end{array}\right].
$$
The shift $\alpha_{V_{14}}$ is $4$.
\end{description}
\end{theorem}

\begin{proof}
By Theorem~\ref{Mukai} these varieties are complete intersections in
$G(2,5)$ (for $V_{10}$) and $G(2,6)$ (for $V_{14}$).

Let $H$ be an effective generator of the Picard group of
grassmannian $G=G(r,n)$. Put $I^{G}=I^G_{H^0}(q)+ I^G_{H^1}(q)\cdot
H \ \, (\mathrm{mod}\ \, H^4(G))$, that is
$$
I^{G}=I^G_{H^0}(q)+ I^G_{H^1}(q)\cdot H + \widetilde{I},
$$
where $\widetilde{I}\in H^{>2}(G)$.

By Corollary~\ref{corollary:free term}
$$
I^{G(2,5)}_{H^0}=1+3q+\frac{19}{32}q^2+\frac{49}{2592}q^3+\frac{139}{884736}q^4+\ldots,
$$
$$
I^{G(2,6)}_{H^0}=1+4q+\frac{3}{4}q^2+\frac{95}{5832}q^3+\frac{865}{11943936}q^4+\ldots.
$$

Theorem~\ref{theorem:HV} enables one to consider the $I$-series of
grassmannians $G(2,5)$ and $G(2,6)$ as series in elementary
symmetric functions $x_1+x_2$ and $x_1x_2$, where $x_1$ and $x_2$
are the Chern roots of the dual to the tautological subbundle. Then
$I^G_{H^1}$ is the coefficient at the linear symmetric function
$x_1+x_2$. Thus,
$$
I^{G(2,5)}_{H^1}=10q+\frac{105}{32}q^2+\frac{3115}{23328}q^3+\frac{6875}{5308416}q^4+\ldots,
$$
$$
I^{G(2,6)}_{H^1}=15q+\frac{609}{128}q^2+\frac{6197}{46656}q^3+\frac{528737}{764411904}q^4+\ldots.
$$

By Corollary~\ref{fano lefschetz}, $\alpha_{V_{10}}=6$ and
$\alpha_{V_{14}}=4$. By formula~(\ref{MF}), taken modulo $H^2$,
$$
I^{V_{10}}=1+10Hq+(39+\frac{67}{2}H)q^2+(220+\frac{3200}{9}H)q^3+
(\frac{6291}{4}+\frac{89387}{48}H)q^4+\ldots \ \, (\mathrm{mod}\ \,
H^2),
$$
$$
I^{V_{14}}=1+5Hq+(16+\frac{31}{4}H)q^2+(2+\frac{1031}{18}H)q^3+
(230+\frac{14863}{96}H)q^4+\ldots \ \, (\mathrm{mod}\ \, H^2).
$$

It is easy to see that the expressions from
Example~\ref{example:expressions} enable one to recover coefficients
of the counting matrix of $Y$ in terms of $I^Y \ \, (\mathrm{mod}\
\, H^2 )$.

\end{proof}

\subsection{Expressions for two-pointed invariants in terms of one-pointed ones}

The method of finding of two- (and more) pointed invariants in terms
of one-pointed ones that we use in Theorem~\ref{third}, may be used
in the general case. More precisely, consider a variety $Y$ of
dimension $n$ and quantum self-dual subalgebra $R\subset H^*(Y)$
with basis $\{\gamma_0,\ldots,\gamma_N\}$. Let $d\in \NN$,
$\langle\tau_i\gamma_j\rangle_k=f_{ij}^k(\langle\gamma_t,
\gamma_r\rangle_s)$ (see Remark~\ref{remark:general_expressions}).
Assume that one-pointed invariants of $Y$ and two-pointed prime ones
that correspond to curves of degree $d$ are known. Then functions
$f_{i,j}^{d+1}$ express one-pointed invariants that correspond to
curves of degree $d+1$ in terms of on prime two-pointed ones that
correspond to curves of degree $d+1$ linearly. Moreover, one can
choose a collection $\{(i,j)\}$ (where $i\in \ZZ_{\geq 0}$, $0\leq j
\leq N$, $i+j=d+n-1$) such that a system of linear equations
$\{f_{ij}^{d+1}(\langle \gamma_k, \gamma_l\rangle_{d+1})=\langle
\tau_i\gamma_j\rangle_{d+1}\}$, is given by nondegenerate
upper-triangular matrix. Thus, by induction on $d$ one can find
polynomial expressions for two-pointed invariants in terms of
one-pointed ones.

\begin{theorem}[\cite{BK}, Theorem $5.2$]
\label{theorem:reconstruction} Let $R\subset H^*(X)$ be the
subalgebra generated by Picard group generator and $I^X=I^X_{R}$.
Then Gromov--Witten invariants of type
$\langle\tau_{i_1}\gamma_1,\ldots, \tau_{i_n}\gamma_n\rangle_d$,
$\gamma_i\in R$ are completely determined by coefficients of
$I$-series $I^X$.
\end{theorem}

In fact, there is a way to recover two-pointed invariants from the
constant term of the variety alone. Consider it for the threefold
case example.

Consider a smooth Fano threefold $Y$ of rank $1$ and index $1$. Put
$I^Y_{H^0}=1+d_2q^2+d_3q^3+d_4q^4+d_5q^5+d_6q^6+\ldots$. Consider a
map $f\colon \Aff^5\rightarrow \Aff^5$, given by polynomials
$f_{i-2,3}^i$ that correspond to invariants $\langle\tau_{i-2}
pt\rangle_i=d_i$ for $i=2,\ldots,6$.
\begin{proposition}
\label{lemma:threefold_equations} The map $f$ is birational.
\end{proposition}

\begin{proof}

Direct computations. Remark that this map is biregular if
$$
-495d_3d_5+261d_2d_3^2-312d_4d_2^2+432d_4^2+56d_2^4\neq 0,
$$
which (a posteriori) holds for smooth Fano threefolds of rank $1$.

\end{proof}

\subsection{Golyshev's conjecture.}

Using the same arguments as in this paper, one can reproduce the
counting matrices for complete intersections in projective spaces
(recovered by Golyshev in~\cite{Go2} from the corresponding $D3$
equations from Givental's Theorem in~\cite{Gi}) and the Fano
threefold $V_5$. The counting matrix for it was found by Beauville
in~\cite{Bea}. Golyshev's conjecture holds for these varieties. It
may also be checked for $V_{12}$ (which is a section of the
orthogonal grassmannian $OG(5,10)$ by a linear subspace of
codimension 7), $V_{16}$ (which is a section of the Lagrangian
grassmannian $LG(3,6)$ by a linear subspace of codimension 3), and
$V_{18}$ (which is a section of the grassmannian of the group $G_2$
by the linear subspace of codimension 2). Quantum multiplication by
the divisor class of such grassmannians is computed in ~\cite{FW}.
Evidence in support of the conjecture for these varieties has been
obtained by Golyshev. Counting matrix for $V_{22}$ was computed by
A.\,Kuznetsov; it also agrees with the prediction. Finally, counting
matrices for double cover of $\PP^3$ branched over a quartic, for
double cover of $\PP^3$ branched over a sextic, and for double cover
of the cone over the Veronese surface branched over a cubic may be
also found by this method. For this use their description as smooth
hypersurfaces in weighted projective spaces and extend Givental's
formula for complete intersections in smooth toric varieties to the
case of complete intersections in singular toric varieties.

\subsection{The case of Picard number strictly greater than 1}
\begin{remark}
    All the methods for varieties with Picard number $1$ described above may be easily generalized to the case of Picard
number greater than $1$. One should use \emph{a multidegree}
$d=(d_1,\ldots,d_r)$ instead of the degree and
\emph{multidimensional} variable $q=(q_1,\ldots,q_r)$ instead of
$q$, and so on.
\end{remark}

The author is grateful to V.\,Golyshev for lots of useful
discussions, to A.\,Gathmann for explanations, to A.\,Givental for
the reference to Hori--Vafa conjecture, and to I.\,Cheltsov,
D.\,Orlov, and C.\,Shramov for comments.

\end{document}